\documentclass[12pt]{article}
\usepackage{amssymb,amsmath}
\usepackage{latexsym}

\newtheorem{thm}{Theorem}[section]
\newtheorem{prop}[thm]{Proposition}
\newtheorem{lem}[thm]{Lemma}
\newtheorem{cor}[thm]{Corollary}

\newtheorem{conj}[thm]{Conjecture}

\begin{document}

\newcommand\tobe[1]{{\checkmark $\bigstar${\small \bf#1}
$\bigstar$}}

\def \F {\bm{F}}
\def \k {\mathbf{k}}
\def \p {\prime}

\def \s {\sigma}
\def \f {\phi}
\def \x {\xi}
\def \l {\lambda}
\def \b {\beta}
\def \bq {\begin{equation}}
\def \eq {\end{equation}}
\def \ba {\begin{array}}
\def \ea {\end{array}}
\def \bt  {\bigtriangleup}
\def \w {\wedge}
\def \cd {\cdot}
\def \sm {\setminus}
\def \a {\alpha}
\def \bc {\begin{center}}
\def \ec {\end{center}}
\def \d {\partial}
\def \ov {\overline}
\title{Commutative Moufang loops and alternative algebras}
\author
{Alexander N.\,Grishkov\thanks{Supported by
 FAPESP and CNPq(Brazil)}\\
 Instituto de Matem\'atica e Estat\'\i stica,\\ Universidade de
S\~ao Paulo, S\~ao Paulo, Brasil\\
e-mail: {\em shuragri@gmail.com}\\ [.5mm]
and\\
Ivan P.\,Shestakov\thanks{Partially supported  by CNPq, grant 304633/03-8 and by FAPESP, grant 05/60337-2.}\\
Instituto de Matem\'atica e Estat\'\i stica,\\ Universidade de
S\~ao Paulo, S\~ao Paulo, Brasil\\
and\\
Sobolev Institute of Mathematics,\  Novosibirsk, Russia\\
e-mail: {\em shestak@ime.usp.br}
}
\maketitle

\begin{abstract}
We compute the orders of free commutative Moufand loops of exponent 3 with
$n\leq 7$ free generators and find  embeddings of such loops into a loop of invertible
elements of the free commutative alternative algebra with identity $x^3=0.$ 
\end{abstract}

\section{Introduction}
\hspace{\parindent}
The variety of Moufang loops is one of the most known and most studied classes of loops.
Among them, the commutative loops are  most investigated.
By definition, a  {\it commutative Moufang loop} (CML) is a loop that satisfies the identities
$$
xy=yx,\quad
xy\cd zx=(xy\cd z)x.
$$

Every CML $M$ satisfies the identity $(x,y,z)^3=1$, where $(x,y,z) = (xy\cdot z) (x\cdot yz)^{-1}$ is the {\it associator} of elements $x,y,z\in M$.
Consequently, $M$ has a normal subloop $M^{(1)}$ of exponent 3 generated by associators such that the quotient loop $M/M^{(1)}$ is an abelian
 group. Besides, the normal subloop $M^3$ of $M$ generated by the set $\{x^3|\,x\in M\}$ lies in the nucleus of $M$, that is, $(M^3,M,M)=1$. This implies that the free $n$-generated CML $\,\tilde L_n\,$ is an extension of the free abelian group ${\bf Z}^n$ by a CML of exponent 3. Moreover, $\tilde L_n^3\cap \tilde L^{(1)}_n=(1)$, which yields that the quotient loop $L_n=\tilde L_n/\tilde L_n^3$ and the original CML $\tilde L_n$ have isomorphic associator subloops $\tilde L_n^{(1)}\simeq L_n^{(1)}$. Therefore, the study of free CML reduces to that of free CML of exponent 3 (CML$_3$).

 \smallskip

 The interest to CML increased after appearing of the book by Yu.Manin  \cite{Ma}, where he related to an arbitrary  cubic form  a corresponding universal CML. Till now, the following fundamental question remains open: whether there exists a cubic form whose universal CML is not an abelian group? In the same book, Manin formulated another question on calculating the order of the free CML$_3$ on $n$ generators \cite[Problem 10.2]{Ma}.

An important approach to solution of this question was done by J.D.H.Smith \cite{Sm3}.
He calculated the order of  $L_n$ provided that the {\it Triple Argument Hypothesis} holds, that is, the normal subloop $TAH(L_n)$ of $L_n$ generated by all associators in $L_n$  with an argument appearing thrice vanishes.  However, we show that the Triple Argument Hypothesis does not hold for $n\geq 7$. More exactly, the associator
$((((x,y,a),a,z),t,v),a,w)\neq 1$ in the free CML$_3$ on free generators $x,y,z,t,v,w,a$.
Note that in \cite{San} it was proved that $TAH(L_n)\neq (1)$ for $n\geq 9$, and the corresponding construction there is more complicated.


\smallskip
Let $|L_n|=3^{\delta (n)}$. The Manin problem asks to calculate $\delta(n)$. The (hypothetic) Smith's formula gives the following values for $\delta(n)$ when $n\leq 7$:
$$
\delta(3)=4,\, \delta(4)=12,\, \delta(5)=49,\, \delta(6)=220,\, \delta(7)=1007.
$$

One of the main results of our paper is the following
\begin{thm}\label{t1} Let $L_n$ be the free CML$_3$ on $n$ generators and $|L_n|=3^{\delta(n)}$. Then
$$
\delta(3)=4,\ \delta(4)=12,\ \delta(5)=49,\ \delta(6)=220,\ \delta(7)= 1014\,\,or\,\,1035,
$$
and $\delta(7)=1014$ if and only if $L_n$ can be embedded into a loop of invertible elements of a unital alternative commutative algebra.\end{thm}

\medskip

Note that the results on the structure of free $CML_3$ can be informally separated in two groups. To the first one we refer the papers, starting with the famous fundamental book of R.Bruck \cite{Br}, where some upper estimates are given for  $\delta(n)$ and for the nilpotency class $k(n)$ of $L_n$.
The central idea of these papers is to deduce corollaries from  the main identity of CML in terms of associators.
The most known and useful result of this type is the one due to R.Bruck which says:
\bq\label{Bruck}
f(x,y,z;,a_1,\ldots,a_n)=((\ldots ((x,y,a_1),a_1,a_2),a_2,\ldots,a_n),a_n,z)
\eq
is symmetric on $a_1,\ldots,a_n$ and skew-symmetric on $x,y,z$.


\medskip

To the second group we refer the papers that give low estimates. Such estimates require exact constructions or examples of CML. Usually, the Grassmann algebra over a field with 3 elements is used for those constructions. For instance, Malbos in \cite{Mal} give the following realization of CML$_3$ on the direct sum $G_1\oplus G_1$ where $G_1$ is the odd part of the Grassmann algebra $G=G_0\oplus G_1$:
\bq\label{gr1}
(x_1,x_2)\cd (y_1,y_2)=(x_1+y_1+x_1y_1(x_2-y_2), x_2+y_2+x_2y_2(y_1-x_1)).
\eq
L.Beneteau  in \cite{Be} used this construction to prove that $\delta(4)=12,\ \delta(5)=49$ and $\delta(6)\geq 214$. However, not every CML$_3$ admits such a realization (see \cite{Sm4}). The importance of the Grassmann algebra  in the theory of $CML_3$ is evidenced by the following problem by J.D.F.Smith
 \cite[Problem 2]{Sm4}:

{\bf Problem.} {\em  Develop a faithful exterior algebra representation for the free commutative Moufang loop
of exponent 3.}

In this paper we give such a representation for the free $CML_3$
with 6 generators (see Theorem \ref{t2} below.)

\smallskip
J.D.F.Smith in papers \cite{Sm1, Sm2} used the following universal construction of CML$_3$:

{\it If $G$ is a group with involution $\s$,\ $M(G)=\{x\in G|\,x^{\s}=x^{-1}\}$ and ${x^3=1}$\  $\forall x\in M(G)$, then $M(G)$ is a CML$_3$ with respect to multiplication $x\cd y=xy^{-1}x$. Moreover, every CML$_3$ can be obtained in this way.}

In fact, this construction is a partial case of {\it groups with triality} \cite{Doro}, and the corresponding group $G$ is called a {\it Fisher group}.
Because of the universal character of this construction, till now nobody could apply it for calculation of $\delta(n)$. In \cite{Sm1, Sm2}, the author used it to prove the exactness of the estimate for nilpotency class of $L_n$ obtained earlier by Bruck, that is, $k(n)=n-1$.
A simpler proof of this equality  was obtained by Malbos \cite{Mal} who used construction
(\ref{gr1}).

\medskip
The present paper can be referred both to the first and to the second group of results on CML.
On the one hand, we obtain an upper estimate for the order of  free  CML$_3$ by associating with it a certain 3-algebra and studying  the associative algebra of right multiplications of this 3-algebra.

On the other hand, we construct a CML$_3$ using the following example of prime alternative commutative algebra from \cite{Sh}:

{\it Let $G=G_0\oplus G_1$ be the Grassmann algebra over a field of characteristic $3$,
 $d:G\rightarrow G$ be the even derivation of $G$; then $G=(G,d)$ becomes an alternative commutative algebra with respect to multiplication:
$$
x\cd y= \left\{
\begin{array}{ll}
	 xy,\  &x\in G_0 \hbox{ or } y\in G_0,\\
	 x^dy-xy^d,\ & x,y\in G_1.
\end{array}\right.
$$}
The loop $U(G,d)$ of invertible elements of $(G,d)$ is a CML which is in fact a CML$_3$. Using this example, we prove that
$$
\delta(5)\geq 44,\ \delta(6)\geq 214,\ \delta(7)\geq 1014.
$$
Note that we obtained these estimates considering only CLM of type $U(A)$ for commutative alternative algebra $A$.
In general, not every Moufang loop can be embedded into a loop $U(A)$ of an alternative algebra $A$ \cite{Sh1}.
For CML, we suggest the following 

\begin{conj}\label{c1}The free CML$_3$ on 7 generators can not be embedded into a loop $U(A)$ for commutative alternative algebra $A$.\end{conj}

More exactly, the following result holds.


\begin{thm}\label{t2} The free CML$_3$ $L_n$ on $n<7$ generators is embedded into a loop $U(A)$ for a unital commutative alternative algebra $A,$  moreover, $L_7$ may be embedded in $U(A)$ if and only if the following identity is true for CML
\bq\label{main-id}\ba{l}
f(a,x,y,z,t;b,c)=\\
((((a,x,y),z,b),t,c),b,c)
((((a,x,z),y,b),t,c),b,c)\\
((((a,x,t),y,b),z,c),b,c)^{-1}
((((a,x,b),y,z),t,c),b,c)\\
((((a,x,c),y,z),t,b),b,c)
((((a,x,b),y,c),z,t),b,c)=1.\ea
\eq
Here $(x,y,z)=(xy)z\cd(x(yz))^{-1}.$
\end{thm}

Note that we proved, in particular, that $|TAH(L_7)|=3^7,$ hence Smith's formula \cite{Sm3} give the correct value of the factor loop
$L_7/TAH(L_7)$  for $n=7$ if and only if the identity (\ref{main-id}) holds in $L_7$.

\section{Constructions}
\hspace{\parindent}
As a main tool in our investigation, we use the following example
of commutative alternative algebra with identity $x^3=0$ over a field $\k$ of characteristic $3$ (see \cite{Sh}).

Let $G$ be the Grassmann algebra over $\k$ with a set of free odd generators
$I=\{i^{(j)}\,|\, i,j=0,1, ...\}$. It is clear that $G$ has a unique derivation $d$ such that
$d(i^{(j)})=i^{(j+1)}$.
We have $G=G_0\oplus G_1$ where  $G_0$ is spanned by $1$ and all products  from $I$
of even length  and  $G_1$ is the linear space spanned by all products from $I$ of odd length.

Denote by $C$ the algebra obtained from $G$ by defining a new multiplication via
$$
x\cd y= \left\{
\begin{array}{ll}
	 xy,\  &x\in G_0 \hbox{ or } y\in G_0,\\
	 x^dy-xy^d,\ & x,y\in G_1.
\end{array}\right.
$$
It was proved in \cite{Sh} that $C$ is a prime alternative algebra. One can easily check that the subalgebra of $C$ generated by $I$ satisfies also the identity $x^3=0$.

Let $C_{n}$ be  the subalgebra of $C$  generated by
$\{0,1, \dots ,n-1\},$ where $i=i^{(0)}.$
\begin{conj}\label{m2}
The algebra $C_n$ is a free commutative alternative algebra with identity $x^3=0$ over
$\k$  on free generators $\{0,1, \dots ,n-1\}.$
\end{conj}

\begin{lem}
If $x,y,z\in C\cap G_1$ then

\bq
\label{2}
(x,y,z)=(x\cd y)\cd z-x\cd(y\cd z)=-d(xyz).
\eq
\end{lem}

{\bf Proof.}
By definition of the multiplication in $C$ we have
\begin{equation}
\label{a1}
 \ba{l}
(x,y,z)=(xd(y)-d(x)y)\cd z-x\cd(yd(z)-d(y)z)=\\
xd(y)z-d(x)yz-xyd(z)+xd(y)z=-d(xyz).
\ea
\end{equation}
$\Box$

Below, we will consider $C$ and $C_n$ as algebras with the trilinear
operation $(x,y,z).$

\smallskip

Let $M$ be a $CML_3$. Define $M^{(0)}=M$, $M^{(n)}=
[M^{(n-1)},M,M]$ where $[x,y,z]=((xy)z) (x(yz))^{-1}.$ Then
$$
K(M)=M/M^{(1)}\oplus M^{(1)}/M^{(2)}\oplus \cdots \oplus M^{(i)}/M^{(i+1)}\oplus \cdots
$$
is an elementary abelian group of exponent 3. It is easy to see that
for $x\in M^{(i)},$ $y\in M^{(j)},$ $z\in M^{(k)},$ we get
$[x,y,z]\in M^{(i+j+k+1)}.$ Hence this group admits the structure
of trilinear algebra over the field ${\bf F}_3$ of 3 elements with
operation induced by $[x,y,z].$

We denote by $K=K(L)$  ($K_n=K(L_n)$) the trilinear algebra that corresponds to the free $CLM_3$ $L$ of countable range (respectively, of rang $n$ with free generators $X=\{x_1,\dots ,x_n\}$).

Finally, let $A_n$ be the free alternative commutative
algebra over  $\k$ with free generators $Y=\{y_1,\dots ,y_n\}$ and identity $x^3=0.$ As above, we will consider $A_n$ as an
algebra with the trilinear operation $(x,y,z).$ By ${\cal A}_n$ (${\cal C}_n$) we denote a 3-subalgebra of 
$A_n$ ($C_n$), generated by $Y$ (\{0,1,\dots ,n-1\}).

\begin{lem}
There exist surjective homomorphisms of trilinear algebras:
$$
\phi_n:K_n\to {\cal A}_n,\,\,\,\psi_n:{\cal A}_n\to {\cal C}_n,
$$
where $\phi_n(x_i)=y_i,$ $\psi_n(y_i)=i-1,$ $i=1,\dots ,n$.
\end{lem}
{\bf Proof.}
Let $X=\{x_1,\ldots,x_n\}$ and $Y=\{y_1,\ldots,y_n\}$ be
the sets of free generators of the free $CML_3$ $L_n$ and the 3-algebra $A_n$
respectively. Consider the algebra $A_n^{\sharp}=A_n\oplus \k\cd 1$ obtained by adjoining the exterior unit element $1$ to $A_n$. Then $A_n^{\sharp}$ is a unital alternative commutative algebra and the set $A_n^*=\{1+a\, |\, a\in A_n\}$ forms a $CML_3$ with respect to multiplication. Let $Q_n$ be the subloop of $A_n^*$ generated by the set $\{1+y_i \,|\,i=1,\ldots, n\}.$
We have a natural epimorphism $f_n:L_n\to Q_n,$ $f_n(x_i)=1+y_i,$
$i=1,\ldots , n,$ which induces an epimorphism
$g_n:K(L_n)\to K(Q_n).$

Consider the direct sum of abelian groups
$$
\overline{A_n}=A_n/A_n^2\oplus A_n^3/A_n^4\oplus \cdots\oplus A_n^{2k-1}/A_n^{2k}\cdots,
$$
We can consider $\overline{A_n}$ as a 3-algebra with multiplication $\{\bar{x},\bar{y},\bar{z}\}=\overline{(x,y,z)}.$
It is easy to see that such defined 3-algebra is isomorphic
to a 3-subalgebra $D$ of $A_n$ which is a sum of all homogeneous components of odd degree with respect to the  free generators $Y:$
$D=\sum_{i=1}\oplus A_n^{(2i-1)}.$ The definition of this isomorphism is obvious: $\lambda_n:\overline{A_n}\to D,$
$\lambda_n(a)=x,$ where $x$ is homogeneous element
such that $\bar{x}=a\, (mod\,A_n^{(2i)}),$ if $x\in A_n^{(2i-1)}.$

Define a homomorphism $\tau_n:K(Q_n)\to \overline{A_n}$ by setting
$\tau_n(1+a)=\bar a_1,$
where $a_1$ is the least term of $a.$
Since in $A_n$ we have
$$
[1+a,1+b,1+c]=(((1+a)(1+b))(1+c))((1+a)((1+b)(1+c))^{-1}= 1+(a,b,c)+d,
$$
where $d$ is a sum of monomials in $a,b,c$ of degree more then 3,
then $\tau_n$ is a homomorphism of 3-algebras.

At last, we define
$\phi_n=g_n\cdot \tau_n\cdot \lambda_n,$
then $\phi_n(K_n)\subseteq D\subset A_n.$

The lemma is proved. $\Box$

Theorems \ref{t1} and \ref{t2} will be obtained as corollaries of the following result.

\begin{thm}\label{ma}
The homomorphisms $\phi_n$ and $\psi_n$ are isomorphisms for $n<7$. Moreover, $\psi_7$ is an isomorphism and
$ker\phi_7$ coincides with the subspace $F(A_7)$ of $A_7$ spanned by all the values of the function $F(a,x,y,z,t;b,c)$ from the left part of (\ref{main-id}).
In particular, $\phi_7$ is an isomorphism if and only if identity (\ref{main-id}) holds in any $CML_3$.
\end{thm}

The idea of the proof is as follows. First, we
calculate dimensions of the algebras $C_n,$ $n<8.$ It is evident that $\dim C_n\leq \dim K_n$.
 We prove that
$\dim K_n\leq \dim C_n$ for $n<7$ and
$\dim K_7\leq \dim C_7$ if (\ref{main-id}) holds in any $CML_3.$

As a result of calculations made in the proof, we get the following

\begin{cor}\label{cor1}
Let $3^{\delta(n)}$ be the order of the free $CML_3$ with $n$
generators. Then
\bq\label{num}
\delta(3)=4,\,\,\delta(4)=12,\,\,\delta(5)=49,\,\,\delta(6)=220.
\eq
Moreover, $\delta(7)=1014$ or
$\delta(7)=1035$   according to whether $(\ref{main-id})$ holds or does not hold  in $CML_3.$
\end{cor}

The Theorem \ref{ma} gives a certain evidence in favor of Conjecture \ref{m2} for $n\leq 7.$

\section{The dimensions of ${\cal C}_n, n<8.$}
\hspace{\parindent}
Let ${\bf a}=(a_1, a_2, ...),$ $a_i\in {\bf Z}_+=\{0,1,2, ...\}$ and
$|{\bf a}|=\sum_iia_i.$
Denote ${\bf Z}^{\omega}=\{{\bf a} \,|\, |{\bf a}|<\infty\}.$
 By $C_n({\bf a})$ we denote the set of 3-words in ${\cal C}$ on $\{0,1,2, ... ,n-1\}$ such that the number
of letters which appear $i$ times in this word is equal to $a_i$. By $c_n({\bf a})$ we denote the dimension of the $\k$-space generated by $C_n({\bf a}).$

 For example, if ${\bf a}=(3,0,0,...)$ then
$C_4({\bf a})=\{(0,1,2),(0,1,3),(0,2,3),(1,2,3)\}$
 and $c_4({\bf a})=4.$ It is clear that
\bq
\label{3}
\dim {\cal C}_{n}=\sum_{{\bf a}}c_n({\bf a}).
\eq

For $n=4$ we have $c_4({\bf a})=0$ if ${\bf a}\notin \{{\bf
a}_1=(1,0,\ldots),\ {\bf a}_2=(3,0,\ldots),\ {\bf a}_3=(3,1,0,\ldots)\}.$
It is clear that $c_4({\bf a}_1)=4,$ $c_4({\bf a}_2)=4.$ Let us calculate the
number $c_4({\bf a}_3).$ Starting with this moment, we will use
the notation $x^{\prime}$ instead of $d(x)$ or $x^d.$ If
$\{i,j,k,p\}=\{0,1,2,3\}$ then by (\ref{2}) we have
$$
v_i=((i,j,k),i,p)=((ijk)^{\prime}ip)^{\prime}=(i^{\prime}jkip)^{\prime}.
$$
Hence $c_4({\bf a}_3)=|\{v_0,v_1,v_2,v_3\}|=4.$
This gives by (\ref{3}):
\begin{prop}
$\dim {\cal C}_4=12.$
\end{prop}

Let ${\bf a}=(a_1,\ldots,a_s,0,\ldots),$ $a_1+\cdots+a_s=m$ and
$P=\{0,1,\ldots,m-1\}=\cup_{i=1}^sP_i,$ $|P_i|=a_i,$ $P_i\cap P_j=\emptyset,$ if $i\not=j.$
We denote $h({\bf a})=\dim V({\bf a}),$ where
$V({\bf a})$ is the space generated by 3-words $v$ such that any
fix letter $i\in P_j$  appear $j$ times in $v$. It is clear that $h({\bf a})$ does not
depend on partitions $P=\cup_{i=1}^sP_i$ but depends only on their types $(a_1,\ldots,a_s)$.

We have an evident formula
\bq\label{ob}
c_n({\bf a})={n\choose{a_1,\ldots,a_s}}h({\bf a}),
\eq
which reduces calculations of $c_n({\bf a})$ to that of  $h({\bf a})$.

\begin{lem}\label{3s}
$h({\bf a})=1$ if ${\bf a}=(3,s,0,\ldots).$
\end{lem}

{\em Proof.} Any 3-word  $v\in V({\bf a})$ has a form:
$$
(xyc_1^{\p}c_1c_2^{\p}c_2\cdots
c_s^{\p}c_sz)^{\p},
$$
where $\{x,y,z\}=\{0,1,2\},\{c_1,\ldots,c_s\}=\{3,4,\ldots,s+2\}$. But
all 3-words of this type are linearly dependent. $\Box$

\smallskip
We will do now some general remarks. 
First, in CML we have the following identity \cite{Br}:
\bq
\label{4}
((x,y,z),a,b)=((x,a,b),y,z)+(x,(y,a,b),z)+(x,y,(z,a,b)).
\eq
With this identity, it is easy to prove (by induction) that any
3-word that contains letter $x$ may be written as a linear
combination of 3-words of the form
$((\ldots((x,a_1,a_2),a_3,a_4),\ldots),a_{2n-1},a_{2n})$ where letter
$x$ appears on the first place.

\smallskip

Second, let $v=((\ldots((x,x_1,x_2),x_3,x_4)\ldots),x_{2n-1},x_{2n})$ be
a 3-word of the form described above. We associate with $v$ an associative word $(x_1x_2)(x_2x_3)\cdots(x_{2n-1}x_{2n})$
on the symbols $x_{2i-1}x_{2i}=-x_{2i}x_{2i-1}.$ Sometimes we will
write it in the form: $x_1x_2.x_2x_3.\cdots .x_{2n-1}x_{2n}.$ We will call this associative word a 
{\em 2-word}.

\begin{prop}\label{h}
(i)\,\,$h({\bf a})=4,$ if ${\bf a}=(5,0,\ldots).$

(ii)\,\,$h({\bf a})=5,$ if ${\bf a}=(5,1,0,\ldots).$

(iii)\,\,$h({\bf a})=6,$ if ${\bf a}=(5,2,0,\ldots).$

(iv)\,\,$h({\bf a})=20,$ if ${\bf a}=(7,0,\ldots).$

(v)\,\,$h({\bf a})=1,$ if ${\bf a}=(6,0,1,0,\ldots).$
\end{prop}

{\em Proof.} We will prove only the last four affirmations. The proof
of the first one is more easy.

1. Let $V$ be a set of 2-words of  type $(5,1)$ on  letters
$a,b,c,d,x$ where only $x$ appear two times. Every 2-word from $V$
may written as $v=ax.xb.cd,$ $w=ax.bc.xd$, or $u=ab.cx.xd.$ As above,
we can suppose that the first letter in the corresponding 3-word
is $o$, that is,
$$
v=((oax)'xb)'cd=o'ax'xbcd+oa'x'xbcd+oax^{(2)}xbcd+oax'xb'cd.
$$
Denote: $(o)=o^{(1)}ax^{(1)}xbcd,$ $(a)=oa^{(1)}x^{(1)}xbcd,$
$(x)=oax^{(2)}xbcd,$ $(b)=oax^{(1)}xb^{(1)}cd,$
$(c)=oax^{(1)}bc^{(1)}d,$ $(d)=oax^{(1)}bcd^{(1)}.$ Then we have

$v=(o)+(a)+(x)+(b),$ $w=-(o)-(a)+(b)+(c)+(x),$ $u=(o)+(a)+(b).$\\
Note that
$v,w,u\in k\{(x),\a_1(o)+\a_2(a)+\a_3(b)+\a_4(c)+\a_5(d) | \a_i\in
k, \a_1+\cdots+\a_5=0\}.$\\
It means  that $\dim C(5,1)=5.$ Here and above we write $C(5,1)$ instead of $C((5,1)).$

Furthermore, every element of type $(5,2)$ is a linear combination of the words
of the following structure:

$v=(((oab)^{(1)}cx)^{(1)}xy)^{(1)}yd,$

$w=(((oax)^{(1)}xy)^{(1)}yb)^{(1)}cd,$

$u=(((oax)^{(1)}bc)^{(1)}xy)^{(1)}yd.$\\
Denote

$(x)=oabcdx^{(2)}xy^{(1)}y,$ $y=oabcdx^{(1)}xy^{(2)}y,$

$(r)=o\cdots r^{(1)}\cdots dx^{(1)}xy^{(1)}y, r\in \{o,a,b,c,d\}.$\\
Then

$v=(o)+(a)+(b),$ $w=(o)+(a)+(b)+(x)+(y),$
$u=-(o)-(a)+(b)+(c)+(x).$\\
Hence
$v,w,u\in k\{(x),(y),\a_1(o)+\cdots+\a_5(d) | \a_1+\cdots+\a_5=0\}$
and $\dim C(5,2)=6.$

The main case ${\bf a}=(7,\ldots),$ denote $h({\bf a})=h(7).$ The
typical element of type 7 has a form: $v=((oab)^{(1)}cd)^{(1)}xy.$
Denote

$(r)=o\cdots r^{(2)}\cdots y,$ $(rt)=o\cdots r^{(1)}\cdots t^{(1)}\cdots y,$$r,t\in
\{o,a,b,c,d,x,y\}.$\\
It is easy to prove that

$ v=(o)+(a)+(b)-(oa)-(ob)+(oc)+(od)-(ab)+
(ac)+(ad)+(bc)+(bd).$\\
Hence $C({\bf a})\subseteq V=V_1\oplus V_2,$ where $V_1=k\{(r)
| r\in \{o,a,b,c,d,x,y\}\}.$ $V_2=k\{(rt) | r,t\in
\{o,a,b,c,d,x,y\}\}.$ $\dim V_1=7, \dim V_2=21,$ $\dim V=28.$
Note that $C({\bf a})/V_2=k\{\a_1(o)+\cdots+\a_7(y) | \a_1+\cdots+\a_7=0\},$
hence $\dim C({\bf a})/V_2=6.$

The space $C({\bf a})\cap V_2$ has a basis of the elements
$w=((oab)^{(1)}c_1d_1)^{(1)}x_1y_1-((oab)^{(1)}c_2d_2)^{(1)}x_2y_2.$
where $\{c_1,d_1,x_1,y_1\}=\{c_2,d_2,x_2,y_2\}=\{c,d,x,y\}.$ We
have

\bq\label{4}\ba{l}
w=(oc_1)+(od_1)-(oc_2)-(od_2)+(ac_1)+(ad_1)+(bc_1)+\\
(bd_1)-(ac_2)-(ad_2)-(bc_2)-(bd_2).
\ea\eq

Let $\phi_i$ be a linear function on $V_2$ such that
$\phi_i(pq)=0,$ if $i\not\in \{p,q\}$ and $\phi_i(pq)=1,$ if $i\in
\{p,q\}.$
Then $w\in W=W_1\cap\cdots \cap W_6,$ where $W_j=\ker\phi_j.$

It is easy to see that $\dim W=14.$ Now it is sufficient to prove
that $V_2\cap C({\bf a})=W.$ If $V_3=\{(pq)\in V_2 | o\not\in \{p,q\}\}$
then by (\ref{4}) we get $w\equiv (oc_1)+(pd_1)-(oc_2)-(od_2) (mod
V_3).$ Hence

$(V_2\cap C({\bf a}))/V_3=k\{\a_1(oa)+\cdots+\a_6(oy) |
\a_1+\cdots+\a_6=0\}$\\
and $\dim V_2\cap C({\bf a})=5+\dim (V_3\cap C({\bf a})).$

Note that $V_3\cap C({\bf a})=W\cap V_3.$ Indeed, $W$ has a basis of
elements of the form: $(ij)+(pq)-(ip)-(jq).$ We will prove that
all elements of this type are contained in $C({\bf a}).$ If in (\ref{4}) we
let $c_2=c_1,$ then the corresponding element has a form:

$u=(od_1)-(od_2)+(ad_1)+(bd_1)-(ad_2)-(bd_2)\in C.$\\
Similarly,

$u^{\p}=(od_1)-(od_2)+(ad_1)+(xd_1)-(ad_2)-(xd_2)\in C,$\\
hence $u-u^{\p}=(bd_1)-(xd_1)-(bd_2)+(xd_2)\in C({\bf a}).$ It means that
$C({\bf a})\cap V_2=W$ and $\dim C({\bf a})=20.$

The last case, ${\bf a}=(6,0,1)).$ All elements of this type are
equal to $\pm v=(((oab)^{(1)}cd)^{(1)}xy)^{(1)}pq.$ Hence,
$\dim C({\bf a})=1.$

\begin{cor}
$(i)\,\,\dim {\cal C}_5=49.$

$(ii)\,\,\dim {\cal C}_6=220.$

$(iii)\,\,\dim {\cal C}_7=1014.$
\end{cor}

{\em Proof.} Let us prove  item (iii).
It is clear that $\dim C_7({\bf a})\not=0$ if and only if
${\bf a}\in \{{\bf a}_i\, |\,i=1,2,\ldots, 11\},$
where ${\bf a}_1=(1,0,\ldots),$ ${\bf a}_2=(3,0,\ldots),$ ${\bf a}_3=(5,0,\ldots),$
${\bf a}_4=(3,1,0,\ldots),$ ${\bf a}_5=(7,0,\ldots),$ ${\bf a}_6=(5,1,0,\ldots),$
${\bf a}_7=(3,2,0,\ldots),$ ${\bf a}_8=(3,3,0,\ldots),$ ${\bf a}_9=(6,0,1,0,\ldots),$
${\bf a}_{10}=(5,2,0,\ldots),$ ${\bf a}_{11}=(3,4,0,\ldots).$

By  Lemma \ref{3s} and  Proposition \ref{h},  we get for
$h_i=h({\bf a}_i):$

$h_1=1,$ $h_2=1,$ $h_3=4,$ $h_4=1,$ $h_5=20,$ $h_6=5,$ $h_7=1,$ $h_8=1,$ $h_9=1,$
$h_{10}=6,$ $h_{11}=1.$

Now by (\ref{ob}) and (\ref{3}) we have 
$$
\dim {\cal C}_7=\sum_{i=1}^{11} \dim C_7({\bf a}_i)=7+35+84+140+20+210+210+140+7+126+35=1014.
$$
$\Box$

\section{The dimensions of $K_n$ for $n<8$.}

Let $k({\bf a})$ be an analogy of $h({\bf a})$ for $K.$ Bruck proved, using (\ref{Bruck}) that

\bq \label{5} k({\bf a})\leq 1,\, if\,{\bf a}=(3,s,0,...) . \eq
If $s=2$ we get by (\ref{Bruck}):
$$
((0,i,x),x,q)=-((0,q,x),x,i).
$$
Linearizing this equality, we have in terms of 2-words:
\bq \label{6} ij.pq+ip.jq+jq.ip+pq.ij=0. \eq
\begin{lem}\label{3-2}
In the 3-algebra $K$ we have the following identity in terms of 2-words:
\bq
\label{7}
ia.ja.pa=0.
\eq
\end{lem}

{\em Proof.}
By (\ref{4}) we have:
$$
A=(((x,i,a),j,a),p,a)=(((x,i,a),p,a),j,a)+((x,i,a),(j,p,a),a)=
$$
$$
-A+((x,i,a),(j,p,a),a),
$$
hence $2A=((x,i,a),(j,p,a),a).$

Note that $A$ is antisymmetric on $x,i,j,p$ hence so is
$((x,i,a),(j,p,a),a)$. Therefore,
$$
((x,i,a),(j,p,a),a)=((j,p,a),(x,i,a),a).
$$
But $((x,i,a),(j,p,a),a)=-((j,p,a),(x,i,a),a)$ since $(x,y,z)$ is
antisymmetric on $x,y,z.$ Then $2A=((x,i,a),(j,p,a),a)=0.$ $\Box$

The following identity is a full linearization of (\ref{7}).
\bq
\label{71}
ia.jb.pc+ia.jc.pb+ib.ja.pc+ib.jc.pa+ic.ja.pb+ic.jb.pa=0. \eq

Now we will estimate the dimensions of $K_n, n<8.$ The idea is as follows.
Let $X$ be a set of associative words on elements $\{ij=-ji \,|\, 1\leq i,j\leq n\}.$ We order $X$ by setting $ij>pq$ for $i<j,p<q$  if $q<j$ or $q=j$ and $i>p$.  Then associative words on $X $ of the same length are lexicographically ordered. We will call a word $v$ {\em irregular} if there exists a word $f\in X$ whose leading term
 is a subword of $v.$ Otherwise we will call this word {\em $X$-regular} or {\em regular}. It is obvious
that the set of all $X$-regular words contains some basis of a
space of words of the same degree modulo the ideal generated by $X$. The following technical lemma
plays a crucial role in the proof of the main result.
\begin{prop}\label{meip}
(i) $k({\bf a})\leq 4,$ if ${\bf a}=(5,0, ...).$

(ii) $k({\bf a})\leq 5,$ if ${\bf a}=(5,1,0, ...).$

(iii) $k({\bf a})\leq 6,$ if ${\bf a}=(5,2,0, ...).$

(iv) $k({\bf a})\leq 20,$ if ${\bf a}=(7,0, ...).$
\end{prop}

{\em Proof.} (i) Let $X$ be a set of words of type (\ref{6}).
The only regular $X-$words on letters $12,13,14,23,24,34$
of degree 2 are the following:
$$
12.34,\,13.24,\,14.23,\,23.14.
$$
Hence $k({\bf a})\leq 4,$ if ${\bf a}=(5,0,...).$

(ii) Now let $Y=X\cup \{f_1,f_2,f_3\}$, where
\bq\label{8}
f_1=iq.ja.pa+ia.jq.pa+ia.ja.pq,\eq
\bq\label{8a}
f_2=pi.ja.aq+pj.ia.aq+qi.ja.ap+qj.ia.ap.\eq
Note that $f_1=0$ is a partial linearization of (\ref{7}), and $f_2=0$ follows from the
identity 
$$
f_3=((x,p,a),a,b),b,q)+((x,q,a),a,b),b,p)=0,
$$
which is a particular case of (\ref{Bruck}).
\begin{lem}\label{3a}
Suppose that $a>4.$ Then the unique $Y-$regular
words of the form $xy.zc.tb,$ $\{x,y,z,t,c,b\}=\{1,2,3,4,a,a\}$ are the
following:

$\{12.3a.4a,\,13.2a.4a,\,14.2a.3a,\,23.1a.4a,\,1a.23.4a\}.$
\end{lem}

{\em Proof.}
Let
$V$ be a space generated by the words $\{12.3a.4a,\,13.2a.4a,\,14.2a.3a,\,23.1a.4a\}.$. We will right $v\equiv w$
if $v-w\in V.$ As $f_1=0$ we have

\bq\label{l1} ia.ja.pq\equiv -ia.jq.pa. \eq Analogously,
\bq\label{l2} ia.ja.pq\equiv -ia.jp.aq. \eq But $ia.ja=-ja.ia,$
hence

\bq\label{l3} ia.jq.pa\equiv -ja.iq.pa. \eq From
(\ref{l1}-\ref{l3}) we can get that all words of the type
$ia.jq.pa$ and $ia.ja.pq$ are contained in $V\oplus
k\{1a.23.4a\}.$
$\Box$

 By Lemma \ref{3a} we get that $k({\bf a})\leq 5,$ if ${\bf a}=(5,1,0,
\dots).$

(iii) Let $Z_1=X\cup \{g_1,g_2\}.$ Where

$g_1=iq.ja.pb+iq.jb.pa+ia.jq.pb+ib.jq.pa+ia.jb.pq+ib.ja.pq$

 is the full
linearizing of $f_1=0$ and

 $g_2=pi.ja.bq+pi.jb.aq+pj.ia.bq+pj.ib.aq+qi.ja.bp+qi.jb.ap+qj.ia.bp+qj.ib.ap$

 is a full linearization of the identity $pi.ia.ap=0,$ the last identity is a corollary
 of (\ref{Bruck}).

Hence a word
$v=ij.pq.rt$ is $Z_1-$irregular if one of the following conditions
hold:

1. One of 2-subwords is $X-$irregular.

2. There exist $s\in \{ij\}$,$l\in \{pq\}$,$m\in \{rt\}$ such that
$s\geq l\geq m.$ In this case we can apply identity $g_1=0.$

3. There exist $\{i_1,j_1,p_1,q_1,r_1,t_1\}$ such that
$\{i_1,j_1\}=\{i= ,j\},$
$\{p_1,q_1\}=\{p,q\},$$\{r_1,t_1\}=\{r,t\},$

moreover, $i_1\geq t_1,$ $j_1\geq p_1,$ $q_1\geq r_1.$
In this case we can apply identity $g_2=0.$

There are $\frac{6!}{8}=90$ words of degree 3 on $\{ij |1\leq i,j\leq 6\}$
of the form $ij.pq.rt,$ where $\{i<j,p<q,r<t\}=\{1,2,3,4,5,6\}.$
We give the list of those words in lexicographical order:
\bq\label{list}\ba{lllll} 1) 12.34.56,&2) 13.24.56,& 3) 23.14.56,&
4) 14.23.56,&
5) 24.13.56,\\
6) 34.12.56,& 7) 12.35.46,& 8) 13.25.46, &
9)23.15.46,&10) 15.23.46,\\
 11) 25.13.46,& 12) 35.12.46,&
13)12.45.36,&14) 14.25.36,&
 15) 24.15.36,\\
  16) 15.24.36&
17)25.14.36,&18) 45.12.36,& 19) 13.45.26,&
 20) 14.35.26,\\
21)34.15.26,&22) 15.34.26,& 23) 35.14.26,& 24) 45.13.26,&
25)23.45.16,\\26) 24.35.16,& 27) 34.25.16,& 28) 25.34.16,&
29)35.24.16,&
30) 45.23.16,\\
 31) 12.36.45,& 32) 13.26.45,&
33)23.16.45,&
34) 16.23.45,& 35) 26.13.45,\\
 36) 36.12.45,&
37)12.46.35,&
38) 14.26.35,& 39) 24.16.35,& 40) 16.24.35,\\
41)26.14.35,&42) 46.12.35,&
 43) 13.46.25,& 44) 14.36.25,&
45)34.16.25,\\
46) 16.34.25,& 47) 36.14.25,&
 48) 46.13.25,&
49)23.46.15,&
50) 24.36.15,\\
 51) 34.26.15,& 52) 26.34.15,&
53)36.24.15,&54) 46.23.15,& 55) 12.26.34,\\
 56) 15.26.34,&
57)25.16.34,&
58) 16.25.34,& 59) 26.15.34,& 60) 56.12.34,\\
61)13.56.24,&62) 15.36.24,&
 63) 35.16.24,& 64) 16.35.24,&
65)36.15.24,\\
66) 56.13.24,& 67) 23.56.14,&
 68) 25.36.14,&
69)35.26.14,&70) 26.35.14,\\ 71) 36.25.14,&
 72) 56.23.14,&
73)14.56.23,&74) 15.46.23,& 75) 45.16.23,\\
 76) 16.45.23,&
77)46.15.23,&
78) 56.14.23,& 79) 24.56.13,& 80) 25.46.13,\\
81)45.26.13,&82) 26.45.13,&
 83) 46.25.13,& 84) 56.24.13,&
85)34.56.12,\\
86) 35.46.12,& 87) 45.36.12,&
 88) 36.45.12,&
89)46.35.12,&90) 56.34.12.\ea\eq
\begin{lem}\label{3-w}
In the list above only $20$ words number
\{1-4,7-10,13-15,19,21,25,31-34,38,56\} are regular.
\end{lem}

{\em Proof.}
 {\bf Step 1}. The words number
$\{5,6,11,12,17,18,23,24,29,30,35-37,41-43,44,47-51,53-55,60-62,64-74,76-90\}$
are not regular since contain some irregular 2-subword.

{\bf Step 2}. The words number $\{16,20,22,26-28,40,46,52,59\}$ are not
regular since have the form $ij.pq.rs, i>p>r.$

{\bf Step 3}. The words number $\{39,45,57,63,75\}$ are not regular since
have the form $ij.pq.rs, i>s,j>p,q>r,i<j,p<q,r<s.$

{\bf Step 4}. The rest $20$ words number
\{1-4,7-10,13-15,19,21,25,31-34,38,56\} are regular. $\Box$

Let $Q$ be a subset of $K$ of all 4-words of the form:
$v=ij.pq.rt.sm,$ where $\{i,j,p,q,r,t,s,m\}=\{1,2,3,4,b,b,c,c\}.$
Clear that $k({\bf a})=\dim k\{Q\}.$

As any 4-word on $\{1,2,3,x,x,y,y,z,z\}$ is
antisymmetric on $\{1,2,3\}$ hence

 \bq\label{s2} xy.xb.yc.bc= 0.
\eq

A linearizing of (\ref{s2}) gives

\bq\label{s3}
g_3=24.1b.3c.bc+23.1b.4c.bc+14.2b.3c.bc+13.2b.4c.bc= 0. \eq

Let $Z_2=Z_1\cup \{g_3\}.$

\begin{lem}\label{l5}
If $Z_2$ as above and we suppose that $c>b>4>3>2>1$ then the set of
$Z_2-$regular words from $Q$ is contained in the following list
\[
\ba{lll}
v_1=12.3b.4c.bc,& v_2=13.2b.4c.bc,& v_3=23.1b.4c.bc,\\
v_4=14.2b.3c.bc,& v_5=1b.23.4c.bc,& v_6=1c.23.4b.bc,\\
v_7=1b.2c.34.bc.& &
\ea
\]
\end{lem}

{\em Proof.}
Let $v=ij.pq.rt.sm\in Q$ and $\{s,m\}\subset \{1,2,3,4\}.$ Then $v$ is not $Z_2$-regular
since the subword $pq.rt.sm$ is not $Z_1$-regular.
If $|\{s,m\}\cap \!\{1,2,3,4\}\!|=1$ then by Lemma \ref{3a} we get that
$v=ij.px.rx.rm,$ $\{x,r\}=\{b,c\}.$
But in this case $v=-ij.px.rm.rx.$ Hence we can suppose that $\{s,m\}=\{b,c\}.$
Then the subword $w=ij.pq.rt$ is $Z_2-$regular, and by Lemma \ref{3-w} we have 20 possibilities
for $w$ numerated in that Lemma as $\{1-4,7-10,13-15,19,21,25,31-34,38,56\}.$

In cases $1-4$ we have $w.bc=0,$ since $bc.bc=0.$

In cases $7,8,9,10,14,34,56,$ we have the words $v_1,v_2,v_3,v_5,v_4,v_6,v_7$ correspondingly.

 In  cases $13,19,25,$ we have the words $-v_1,-v_2,-v_3$ correspondingly. For example,
 (the case 13): $12.4b.3c.bc=-12.3b.4c.bc=-v_1.$

 In cases $15,21,$ we can apply the relation $g_3=0.$

 At last, in  cases $31,32,33,56,$ we can apply  the relation $ic.jb.bc=-ib.jc.bc$ to a certain subword\textsl{}.
 Lemma \ref{l5} is proved.
 $\Box$

 By Lemma \ref{l5}, we get that $k(5,2)\leq 7.$ Therefore, Proposition \ref{meip} is proved. $\Box$
 
 We can now obtain Theorem \ref{t1} as a corollary of Propositions \ref{h} and \ref{meip}.
 
 {\bf Proof} of Theorems \ref{ma} and \ref{t2}. We proved in  Lemma \ref{l5} that $k(5,2)\leq 6$ if the 2-words
 $\{v_1,\dots,v_7\}$ are linearly dependent in $K_7$ which is the same that in $L_7.$
 But in $C_7$ we have a unique linear relation between the elements $\{v_1,\ldots,v_7\}$,
 namely $v_1+v_2-v_4+v_5+v_6+v_7=0.$ It is clear that this relation is equivalent to identity
 (\ref{main-id}), moreover, the identity (\ref{main-id}) holds in $U(A)$ by the following lemma.
\begin{lem}\label{Alb}
Let $Alt_7$ be the free commutative alternative algebra with identity $x^3=0.$
Then in $Alt_7$ we have the following identity (in notation of Lemma \ref{l5})
\bq\label{to2}
v_1+v_2-v_4+v_5+v_6+v_7=0.
\eq
\end{lem}
{\bf Proof.} For the proof we used the programm
ALBERT \cite{Albert} for calculations in nonassociative algebras.
 $\Box$
 
Hence $\psi_7$ is an isomorphism and Theorem \ref{ma} is proved.
At last, Theorem \ref{t2} hence from Theorem \ref{ma}.

\end{document}